\documentclass[12pt,leqno]{article}
\usepackage{amsmath,amssymb,amscd,latexsym}

\usepackage[all]{xy}
\newcommand{\C}{{\mathbb{C}}}
\newcommand{\Q}{{\mathbb{Q}}}
\newcommand{\R}{{\mathbb{R}}}
\newcommand{\Z}{{\mathbb{Z}}}
\newcommand{\cover}{\mathrm{cover}}
\newcommand{\ddet}{\mathrm{det}}
\newcommand{\of}{\overline{f}}
\newcommand{\sep}{\mathrm{sep}}
\newcommand{\spa}{\mathrm{span}\,}
\newcommand{\tr}{\mathrm{tr}}
\newcommand{\End}{\mathrm{End}\,}
\newcommand{\Ker}{\mathrm{Ker}\,}
\newcommand{\rtop}{\mathrm{top}\,}
\newcommand{\Ah}{{\mathcal A}}
\newcommand{\Bh}{\mathcal{B}}
\newcommand{\Fh}{{\mathcal F}}
\newcommand{\Gh}{{\mathcal G}}
\newcommand{\Nh}{{\mathcal N}}
\newcommand{\Ph}{\mathcal{P}}
\newcommand{\Qh}{\mathcal{Q}}
\newcommand{\ea}{\mathfrak{a}}
\newcommand{\eU}{\mathfrak{U}}
\newcommand{\oa}{\overline{a}}
\newcommand{\oz}{\overline{z}}
\newcommand{\hA}{\hat{A}}
\newcommand{\tx}{\tilde{x}}
\newcommand{\ty}{\tilde{y}}
\newcommand{\tY}{\tilde{Y}}
\newcommand{\ohne}{\setminus}
\newcommand{\silo}{\stackrel{\sim}{\longrightarrow}}
\newcommand{\tei}{\, | \,}

\newcommand{\halb}{\frac{1}{2}}
\newtheorem{theorem}{Theorem}

\newtheorem{prop}[theorem]{Proposition}
\newtheorem{cor}[theorem]{Corollary}

\newenvironment{rem}{\noindent {\bf Remark}}{}
\newenvironment{rems}{\noindent {\bf Remarks}}{}
\newenvironment{theoremon}{\noindent {\bf Theorem}\it}{}

\newenvironment{proof}{\noindent {\bf Proof}}{\mbox{}\hspace*{\fill}$\Box$}
\begin{document}
\title{Fuglede--Kadison determinants and entropy for actions of discrete amenable groups}
\author{Christopher Deninger}
\date{\today}
\maketitle

\section{Introduction}
\label{sec:1}

Consider a discrete group $\Gamma$ and an element $f$ in the integral group ring $\Z \Gamma$. Then $\Gamma$ acts from the left on the discrete additive group $\Z \Gamma / \Z \Gamma f$ by automorphisms of groups. Dualizing, we obtain a left $\Gamma$-action on the compact Pontrjagin dual group $X_f = \widehat{\Z \Gamma / \Z \Gamma f}$ by continuous automorphisms of groups. By definition, $X_f$ is a closed subshift of $(S^1)^{\Gamma}$.

For measurable or topological actions of amenable groups a good entropy theory is available, \cite{OW}, \cite{M}. Thus from now on let $\Gamma$ be a finitely generated, discrete amenable group. We are interested in determining the topological entropy $h_f$ of the $\Gamma$-action on $X_f$. It is known that $h_f$ agrees with the measure theoretic entropy of the $\Gamma$-action on $X_f$ equipped with its Haar probability measure.

In their paper \cite{FK} from 1952, Fuglede and Kadison introduced a determinant for units in finite factor von Neumann algebras. More recently, L\"uck has generalized their theory to non-units and to group von Neumann algebras which are not necessarily factors \cite{L}. He needed this to define $L^2$-torsion for general covering spaces.

The von Neumann algebra $\Nh\Gamma$ of the amenable group $\Gamma$ is a certain completion of $\C\Gamma$ and hence we may consider the Fuglede--Kadison--L\"uck determinant $\det_{\Nh\Gamma} f$ of $f$.

We refer to section \ref{sec:6} for the definition of a strong F{\o}lner sequence. It is known that finitely generated virtually nilpotent groups have strong F{\o}lner sequences.

Our main result is this:

\begin{theorem}
\label{t1_neu}
  The formula \quad $h_f = \log \ddet_{\Nh\Gamma} f$ \quad holds under the following assumptions on $f$ and $\Gamma$:\\
a) The group $\Gamma$ has a strong F{\o}lner sequence.\\
b) $f$ is a unit in the $L^1$-convolution algebra $L^1 (\Gamma)$\\
c) $f$ is positive in $\Nh\Gamma$.
\end{theorem}

Condition a) is not necessary for the inequality $h_f \ge \log \det_{\Nh\Gamma} f$. It occurs because at one point in the proof of the inequality $h_f \le \log \det_{\Gamma} f$ we need to transfer an estimate which is only available for $L^2$-norms to an estimate of $L^{\infty}$-norms. The $L^2$-context is natural for von Neumann algebras whereas the $L^{\infty}$-context is closer to the definition of entropy (in terms of $\varepsilon$-separated sets). \\
Condition b) is in general stronger than invertibility of $f$ in $\Nh\Gamma$. However for certain classes of groups and in particular for virtually nilpotent groups, invertibility in $L^1 (\Gamma)$ and $\Nh\Gamma$ are equivalent. Dynamically, condition b) implies that the $\Gamma$-action on $X_f$ is expansive.\\
Finally, as explained in section \ref{sec:7}, condition c) can be removed if $h_f = h_{f^*}$ and if a Yuzvinskii addition formula holds for the entropies of $\Gamma$-actions. Here $^*$ is the canonical involution on $\Z\Gamma$. For $\Gamma \cong \Z^n$ this formula is known and it is expected to hold in general. 

Note that many elements $f$ satisfying conditions b) and c) can be constructed as follows: Choose $g$ in $\Q \Gamma$ with $\| g \|_1 < 1$. Then $1 + g$ is a unit in $L^1 (\Gamma)$ and for a suitable integer $N \ge 1$ the element $h = N (1+g)$ lies in $\Z \Gamma \cap L^1 (\Gamma)^{\times}$. The element $f = h h^*$ then satisfies both b) and c).

Let us now look at the abelian case $\Gamma \cong \Z^n$ which motivated the paper. In this case, a) is satisfied and c) can be dropped. Moreover, by theorems of Schmidt and Wiener condition b) is equivalent to the $\Z^n$-action on $X_f$ being expansive. Viewing $f$ in $\Z [\Z^n]$ as a Laurent polynomial in $n$ variables we have according to \cite{L} Example 3.13:
\begin{equation}
  \label{eq:1}
  \log \ddet_{\Nh\Gamma} f = \int_{T^n} \log |f (z)| \, d\mu (z) \; .
\end{equation}
Here $T^n$ is the real $n$-torus and $\mu$ its normalized Haar measure. The integral on the right is called the (logarithmic) Mahler measure $m (f)$ of $f$. Thus theorem \ref{t1_neu} becomes the equality 
\begin{equation}
  \label{eq:2}
  h_f = m (f)
\end{equation}
if $X_f$ is expansive or equivalently if $f$ does not vanish in any point of $T^n$. Formula (\ref{eq:2}) was in fact proved for {\it all} $f \neq 0$ by Lind, Schmidt and Ward using a completely different method in \cite{LSW}. Their method has the advantage of working in the non-expansive case as well but it cannot be generalized to non-abelian groups $\Gamma$.

The first interesting non-abelian group to which theorem \ref{t1_neu} applies is the $3$-dimensional integral Heisenberg group. In this case $\Nh\Gamma$ is of type $\mathrm{II}_1$. 

Our proof of theorem \ref{t1_neu} is a kind of infinite dimensional generalization of the argument for finite groups $\Gamma$ which can be found in section \ref{sec:7}. For the generalization we use ideas from Rita Solomyak's beautiful paper \cite{S}. There she explains the occurrence of the Mahler measure as the common entropy of two different types of dynamical systems by using the theory of electrical networks. All polynomials in \cite{S} vanish on $T^n$, so strictly speaking the intersection of both papers is empty. However the spirit of section 3 of her paper pervades our approach.

Another important ingredient in the proof of our entropy formula is an approximation lemma for traces on group von Neumann algebras that was proved by L\"uck in his solution of a problem of Gromov on $L^2$-Betti numbers, \cite{L} Lemma 13.42.

The Mahler measure $m (f)$ appears in several branches of mathematics, c.f. \cite{B1} and \cite{B2} and their references. In particular $m (f)$ is interesting from an arithmetic point of view. If $f$ does not vanish on $T^n$ for instance, I showed in \cite{D} that $m (f)$ is a Deligne period of a certain mixed motive over $\Q$. One may wonder whether some non-commutative analogue of this fact holds for the values $\log \det_{\Nh\Gamma} f$ if $f$ is in $\Z \Gamma \cap L^1 (\Gamma)^{\times}$.

I would like to thank Wolfgang L\"uck for presenting me with a copy of his book on $L^2$-invariants. This started the project. I am also grateful to Siegfried Echterhoff for numerous useful discussions and to Viktor Losert for a helpful letter on some aspects of harmonic analysis. Finally, I would like to thank the Max-Plack-Institute in Bonn  for its support where some of the work on this paper was done.

\section{Background on entropy}
\label{sec:2}

In this section we briefly review known facts about entropy of $\Gamma$-actions which are well documented for $\Z^n$-actions e.g. in \cite{Sch} but for which I could find only partial references for general $\Gamma$. The different notions of entropy of topological, metric and measure theoretic nature coincide for the dynamical systems we are interested in. We also give a convenient description of entropy for subshifts. 

As before, let $\Gamma$ be a finitely generated discrete group which is amenable c.f. \cite{P} or the survey \cite{L} 6.4.1. This is equivalent to the existence of a ``(right) F{\o}lner sequence'' $F_1 , F_2 , \ldots$ of finite subsets of $\Gamma$ such that for every finite subset $K$ of $\Gamma$ we have
\[
\lim_{n\to \infty} \frac{|F_n K \bigtriangleup F_n|}{|F_n|} = 0 \; .
\]
Now assume that $\Gamma$ operates from the left by homeomorphisms on a compact metric space $(X,d)$. For $F \subset \Gamma$, a subset $E$ of $X$ is called $(F, \varepsilon)$-separated if for all $x \neq y$ in $E$ there exists $\gamma \in F$ with $d (\gamma^{-1} x , \gamma^{-1} y) \ge \varepsilon$. Let $s_F (\varepsilon)$ be the maximum of the cardinalities of all $(F, \varepsilon)$-separated subsets.

A subset $E$ of $X$ is called $(F, \varepsilon)$-spanning if for every $x$ in $X$ there is some $y$ in $E$ such that we have $d (\gamma^{-1} x , \gamma^{-1} y) < \varepsilon$ for all $\gamma$ in $F$. Let $r_F (\varepsilon)$ be the minimum of the cardinalities of all $(F, \varepsilon)$-spanning sets. It is finite and because of the inequalities, \cite{E},~\S\,3
\[
r_F (\varepsilon) \le s_F (\varepsilon) \le r_F (\varepsilon/2)
\]
the quantity $s_F (\varepsilon)$ is finite as well.

It follows that
\[
h_{\sep} := \lim_{\varepsilon\to 0} \varlimsup_{n\to \infty} \frac{1}{|F_n|} \log s_{F_n} (\varepsilon)
\]
and
\[
h_{\spa} := \lim_{\varepsilon\to 0} \varlimsup_{n\to \infty} \frac{1}{|F_n|} \log r_{F_n} (\varepsilon)
\]
agree. Note that the $\varepsilon$-limits exist by monotonicity.

For an open cover $\eU$ of $X$ let $N (\eU)$ be the cardinality of a minimal subcover of $\eU$. For a finite subset $F$ of $\Gamma$ let
\[
\eU^F = \bigvee_{\gamma\in F} \gamma \eU
\]
be the common refinement of the finitely many covers $\gamma \eU$. Using \cite{LW1} Theorem 6.1 one sees that the limit
\[
h (\eU) = \lim_{n\to\infty} \frac{1}{|F_n|} \log N (\eU^{F_n})
\]
exists and is independent of the F{\o}lner sequence $(F_n)$. The quantity $h_{\cover} = \sup_{\eU} h (\eU)$ is the topological entropy of the $\Gamma$-action on $X$ introduced in \cite{M} Ch. 5. 

Let $\eU_{\varepsilon}$ be the open cover of $X$ consisting of all open $\varepsilon$-balls. As in the proof of proposition 13.1 of \cite{Sch} one shows the inequalities 
\[
N (\eU^F_{\varepsilon}) \le r_F (\varepsilon) \le s_F (\varepsilon) \le N (\eU^F_{\varepsilon/2}) \; .
\]
They imply the following standard result:

\begin{prop}
  \label{t21}
  \begin{eqnarray*}
    h_{\rtop} := h_{\cover} = h_{\sep} = h_{\spa} & = & \lim_{\varepsilon\to 0} \varliminf_{n\to \infty} \frac{1}{|F_n|} \log r_{F_n} (\varepsilon) \\
 & = & \lim_{\varepsilon\to 0} \varliminf_{n\to\infty} \frac{1}{|F_n|} \log s_{F_n} (\varepsilon) \; .
  \end{eqnarray*}
\end{prop}

The notion of the measure theoretic entropy for actions of amenable groups was introduced in \cite{OW}. Assume our $\Gamma$ operates from the left by measure preserving automorphisms on a probability space $(X, \Ah , \mu)$. For a finite measurable partition $\Ph$ of $X$ and finite $F$ in $\Gamma$ define $\Ph^F$ as above and set
\[
h (\Ph) = \lim_{n\to\infty} \frac{1}{|F_n|} H (\Ph^{F_n})
\]
where
\[
H (\Qh) = - \sum_{A \in \Qh} \mu (A) \log \mu (A) \; .
\]
Again the limit exists by (a simpler version of) \cite{OW} Theorem 6.1 and it is independent of the choice of $(F_n)$. Set $h_{\mu} = \sup_{\Ph} h (\Ph)$. 

\begin{theorem}
  \label{t22}
Let $X$ be a compact metrizable abelian group with (normalized) Haar-measure $\mu$ on which $\Gamma$ acts by automorphisms. Then we have $h_{\rtop} = h_{\mu}$. 
\end{theorem}

\begin{proof}
  Choose an invariant metric $d$ on $X$ inducing the topology and let $E$ be an $(F,\varepsilon)$-separated subset. Define the metric $d_F$ on $X$ by $d_F (x,y) = \max_{\gamma \in\Gamma} d (\gamma^{-1} x , \gamma^{-1}y)$ and let $U = U_{d_F} (0,\varepsilon)$ be the open ball of radius $\varepsilon$ in the $d_F$-metric around $0 \in X$. Then the sets $x + U$ for $x \in E$ are pairwise disjoint. Hence we have
\[
|E| \mu (U) = \mu \left( \Dot{\bigcup_{x \in E}} x + U \right) \le \mu (X) = 1 \; ,
\]
and therefore
\begin{equation}
  \label{eq:3}
  s_F (\varepsilon) = \max |E| \le \mu (U)^{-1} \; .
\end{equation}
Let $\Ph$ be a finite partition of $X$ into measurable sets of $d$-diameter $< \varepsilon$. Then every set in $\Ph^F$ has $d_F$-diameter $< \varepsilon$ and hence is contained in $x + U$ for some $x$ in $X$. This gives:
\[
H (\Ph^F) = - \sum_{A \in \Ph^F} \mu (A) \log \mu (A) \ge - \sum_{A \in \Ph^F} \mu (A) \log \mu (U) \overset{(\ref{eq:3})}{\ge} \log s_F (\varepsilon) \; .
\]
Hence we have
\[
h (\Ph) \ge \varlimsup_{n\to\infty} \frac{1}{|F_n|} \log s_{F_n} (\varepsilon) \; .
\]
For every $\varepsilon > 0$ there is a finite measurable partition into sets of diameter $<\varepsilon$. This implies that
\[
h_{\mu} = \sup_{\Ph} h (\Ph) \ge \sup_{\varepsilon >0} \varlimsup_{n\to\infty} \frac{1}{|F_n|} \log s_{F_n} (\varepsilon) = h_{\sep} = h_{\rtop} \; .
\]
On the other hand, according to \cite{M}, 5.1.3 we have $h_{\rtop} \ge h_{\nu}$ for {\it all} $\Gamma$-invariant Radon probability measures $\nu$ on $X$.
\end{proof}

For closed subshifts we will now describe a more convenient way to calculate $h_{\sep}$ and $h_{\spa}$. This will be used in all our later calculations.

Let $G$ be a compact topological group with a left $G$-invariant metric $\vartheta$. The amenable group $\Gamma$ acts from the left on $G^{\Gamma}$ by the formula $(\gamma x)_{\gamma'} = x_{\gamma^{-1} \gamma'}$ for any $x = (x_{\gamma'})$ in $G^{\Gamma}$. If we view $G^{\Gamma}$ as the group (under componentwise multiplication) of formal series $G [|\Gamma|]$ then the action by $\Gamma$ is left multiplication:
\[
\gamma \sum_{\gamma'} x_{\gamma'} \gamma' = \sum_{\gamma'} x_{\gamma'} \gamma \gamma' = \sum_{\gamma'} x_{\gamma^{-1} \gamma'} \gamma' \; .
\]
Let $X$ be a closed $\Gamma$-invariant subspace of $G^{\Gamma}$. For $F \subset \Gamma$ a subset $E$ of $X$ is called $[F,\varepsilon]$-separated if for all $x \neq y$ in $E$ there is some $\gamma$ in $F$ with $\vartheta (x_{\gamma} , y_{\gamma}) \ge \varepsilon$. It is called $[F,\varepsilon]$-spanning if for every $x$ in $X$ there exists a $y$ in $E$ with $\vartheta (x_{\gamma} , y_{\gamma}) < \varepsilon$ for all $\gamma$ in $F$.

Let $s'_F (\varepsilon)$ resp.~$r'_F (\varepsilon)$ be the maximum resp.~minimum of the cardinalities of all $[F , \varepsilon]$-separated resp.~$[F,\varepsilon]$-spanning subsets $E$ of $X$.

Fix a sequence $(c_{\gamma})$ of positive real numbers indexed by $\gamma \in \Gamma$ such that $c_e = 1$ and $\sum c_{\gamma} < \infty$. Then the formula
\[
d (x,y) = \sum_{\gamma \in \Gamma} c_{\gamma} \vartheta (x_{\gamma} , y_{\gamma})
\]
defines a left-invariant metric on $G^{\Gamma}$ which induces the topology of $G^{\Gamma}$. The inequality $d (x,y) \ge \vartheta (x_e , y_e)$ implies that $[F,\varepsilon]$-separated sets are $(F,\varepsilon)$-separated (for $d$) and $(F,\varepsilon)$-spanning sets are $[F,\varepsilon]$-spanning. It follows that we have $s'_F (\varepsilon) \le s_F (\varepsilon)$ and $r'_F (\varepsilon) \le r_F (\varepsilon)$. The next result generalizes \cite{Sch} Proposition 13.7.

\begin{prop}
  \label{t23}
For a closed $\Gamma$-invariant subshift $X \subset G^{\Gamma}$ as above we have the following formulas for the topological entropy
\begin{eqnarray*}
  h_{\rtop} & = & \lim_{\varepsilon\to 0} \varliminf_{n\to\infty} \frac{1}{|F_n|} \log r'_{F_n} (\varepsilon) = \lim_{\varepsilon\to 0} \varlimsup_{n\to\infty} \frac{1}{|F_n|} \log r'_{F_n} (\varepsilon) \\
 & = & \lim_{\varepsilon\to 0} \varliminf_{n\to\infty} \frac{1}{|F_n|} \log s'_{F_n} (\varepsilon) = \lim_{\varepsilon\to 0} \varlimsup_{n\to\infty} \frac{1}{|F_n|} \log s'_{F_n} (\varepsilon) \; .
\end{eqnarray*}
\end{prop}

\begin{proof}
  \quad {\bf 1} For every $\varepsilon > 0$ there exist $\varepsilon' > 0$ and a finite subset $K \subset \Gamma$ such that for all $x,y$ in $X$ with $d (x,y) \ge \varepsilon'$ there is some $\gamma$ in $K$ with $\vartheta (x_{\gamma} , y_{\gamma}) \ge \varepsilon$. 

Namely set $A = \sum c_{\gamma}$ and choose a finite subset $K \subset \Gamma$ with $\sum_{\gamma \in \Gamma \ohne K} c_{\gamma} < \varepsilon$. Set $\varepsilon' = (A+D) \varepsilon$ where $D$ is the diameter of $G$ in the $\vartheta$-metric. If we have $\vartheta (x_{\gamma} , y_{\gamma}) < \varepsilon$ for all $\gamma \in K$ then
\[
d (x,y) < \sum_{\gamma \in K} c_{\gamma} \varepsilon + D \sum_{\gamma \in \Gamma \ohne K} c_{\gamma} \le A \varepsilon + D \varepsilon = \varepsilon' \; .
\]

{\bf 2} Fix $\varepsilon > 0$ and choose $\varepsilon' , K$ as in {\bf 1}. For any finite subset $F \subset \Gamma$ set $F' = \{ f \in F \tei fK \subset F \}$. Then any $(F' , \varepsilon')$-separated set $E$ in $X$ is $[F , \varepsilon]$-separated.

Namely, consider elements $x \neq y$ in $E$. By assumption there is some $f$ in $F'$ with $d (f^{-1} x , f^{-1} y) \ge \varepsilon'$. By {\bf 1} there is some $\gamma$ in $K$ with $\vartheta (x_{f\gamma} , y_{f\gamma}) = \vartheta ((f^{-1} x)_{\gamma} , (f^{-1} y)_{\gamma}) \ge \varepsilon$. Since $f\gamma \in F$ it follows that $E$ is $[F,\varepsilon]$-separated.

{\bf 3} Let $(F_n)$ be a F{\o}lner sequence. Then for fixed $K$ the sequence $(F'_n)$ with $F'_n = \{ f \in F_n \tei fK \subset F_n \}$ is a F{\o}lner sequence as well.

Namely, the F{\o}lner property for a sequence $(F_n)$ is equivalent to
\[
\lim_{n \to\infty} \frac{|F_n \gamma \ohne F_n|}{|F_n|} = 0 \quad \mbox{for all} \; \gamma \in \Gamma \; .
\]
This property is readily checked for $F'_n$ since we have
\[
|F'_n \gamma \ohne F'_n| \le |F_n \gamma \ohne F'_n| \le |F_n \gamma \ohne F_n| + |F_n \ohne F'_n|
\]
and since $\lim_{n\to\infty} |F'_n| / |F_n| = 1$ because $(F_n)$ was F{\o}lner.

{\bf 4} Given $\varepsilon > 0$ choose $\varepsilon' , K$ as in {\bf 1} and consider $(F_n)$ and $(F'_n)$ as in {\bf 3}. Using {\bf 2} it follows that we have
\[
s_{F'_n} (\varepsilon') \le s'_{F_n} (\varepsilon)
\]
and hence
\[
\varliminf_{n\to\infty} \frac{1}{|F'_n|} \log s_{F'_n} (\varepsilon') \le \varliminf_{n\to \infty} \frac{1}{|F'_n|} \log s'_{F_n} (\varepsilon) = \varliminf_{n\to\infty} \frac{1}{|F_n|} \log s'_{F_n} (\varepsilon) \; .
\]
For $\varepsilon \to 0$ the $\varepsilon' = (A+D) \varepsilon$ tends to zero as well. Since the limits exist by monotonicity this shows using proposition \ref{t21} that
\[
h_{\rtop} \le \lim_{\varepsilon\to 0} \varliminf_{n\to \infty} \frac{1}{|F_n|} \log s'_{F_n} (\varepsilon) \; .
\]
The reverse inequality follows from the inequality $s'_{F_n} (\varepsilon) \le s_{F_n} (\varepsilon)$ noted above. Similar arguments for $\varlimsup$ and $r ,r'$ conclude the proof.
\end{proof}
%\newpage
%\input{sec3_fuglede}

\section{Group von Neumann algebras and Fuglede--Kadison determinant}
\label{sec:3}

In this section we first recall the definitions of group von Neumann algebras and the Fuglede--Kadison determinant following \cite{L}. In the amenable case we then prove a formula for this determinant as a limit of finite dimensional determinants. 

For a discrete group $\Gamma$ let $L^2 (\Gamma)$ be the Hilbert space of square summable complex valued functions $x : \Gamma \to \C$. The group $\Gamma$ acts isometrically from the left by the formula $(\gamma x)_{\gamma'} = x_{\gamma^{-1} \gamma'}$. If we represent the elements of $L^2 (\Gamma)$ as formal sums $\sum x_{\gamma'} \gamma'$ this corresponds to left multiplication by $\gamma$. For a Banach space $H$ let $\Bh (H)$ be the algebra of bounded linear operators from $H$ to itself. The group von Neumann algebra $\Nh\Gamma$ of $\Gamma$ is the algebra of $\Gamma$-equivariant bounded operators from $L^2 (\Gamma)$ into itself:
\[
\Nh\Gamma = \Bh (L^2 (\Gamma))^{\Gamma} \; .
\]
The von Neumann trace on $\Nh\Gamma$ is the linear form
\[
\tr_{\Nh\Gamma} : \Nh\Gamma \longrightarrow \C
\]
defined by $\tr_{\Nh\Gamma} (g) = (g (e) , e)$ where $e \in \Gamma \subset L^2 (\Gamma)$ is the unit of $\Gamma$.

The group $\Gamma$ acts isometrically on $L^2 (\Gamma)$ from the right by $(x\gamma)_{\gamma'} = x_{\gamma' \gamma^{-1}}$. This corresponds to right multiplication by $\gamma$ if we view elements of $L^2 (\Gamma)$ as formal sums. Define the operator $R_{\gamma} : L^2 (\Gamma) \to L^2 (\Gamma)$ by $R_{\gamma} (x) = x\gamma$. It is $\Gamma$-equivariant and hence defines an element of $\Nh\Gamma$. The $\C$-algebra homomorphism
\begin{equation}
  \label{eq:4}
  r : \C \Gamma \longrightarrow \Nh\Gamma \; , \; \sum a_{\gamma} \gamma \longmapsto \sum a_{\gamma} R_{\gamma^{-1}}
\end{equation}
is injective and will often be viewed as an inclusion in the following. Restricted to $\C\Gamma$ the von Neumann trace satiesfies $\tr_{\Nh\Gamma} (\sum a_{\gamma} \gamma) = a_e$. Since $R^*_{\gamma} = R_{\gamma^{-1}}$ we have $r (f)^* = r (f^*)$ where for $f = \sum a_{\gamma} \gamma$ we set $f^* = \sum \oa_{\gamma} \gamma^{-1}$.

For a unit $u$ in $\Nh\Gamma$ the element $uu^*$ is positive and applying the functional calculus in $\Bh (L^2 (\Gamma))$ we obtain $\log uu^*$ in $\Nh\Gamma$. The Fuglede--Kadison--L\"uck determinant of $u$ is the positive real number 
\begin{equation}
  \label{eq:5}
  \ddet_{\Nh\Gamma} u = \exp \left( \halb \tr_{\Nh\Gamma} (\log uu^*) \right) \; .
\end{equation}
This defines a homomorphism
\[
\ddet_{\Nh\Gamma} : (\Nh\Gamma)^{\times} \longrightarrow \R^*_+
\]
with $\det_{\Nh\Gamma} (u) = \det_{\Nh\Gamma} (u^*)$. If $u$ is positive, then $u = \sqrt{uu^*}$ and hence we have: 
\begin{equation}
  \label{eq:6}
  \ddet_{\Nh\Gamma} (u) = \exp (\tr_{\Nh\Gamma} (\log u)) \; .
\end{equation}
For an element $f$ of $\C\Gamma$ we also write $\det_{\Nh\Gamma} f$ for $\det_{\Nh\Gamma} r (f)$.

The product in $\C\Gamma$ extends by continuity to its completion in the $\| \; \|_1$-norm and one obtains the $L^1$-convolution algebra $L^1 (\Gamma)$ of $\Gamma$. Right multiplication on $L^p (\Gamma)$ is continuous for $p \ge 1$ including $p = \infty$ because of the estimate $\|\varphi \cdot f \|_p \le \|f \|_1 \, \|\varphi\|_p$ for all $\varphi$ in $L^p (\Gamma)$ and $f$ in $L^1 (\Gamma)$. 

Similarly as before, one obtains injective algebra homomorphisms
\begin{equation}
  \label{eq:7}
  r^{(p)} : L^1 (\Gamma) \longrightarrow \Bh (L^p (\Gamma))^{\Gamma} \quad \mbox{with} \; \|r^{(p)} (f)\| \le \|f \|_1 \; ,
\end{equation}
by mapping $f$ to right multiplication with $f^*$. Note that if we view the elements of our $L^p$-spaces not as formal sums but as functions then multiplication becomes convolution. For $p = 2$ we obtain an injection:
\begin{equation}
  \label{eq:8}
  r = r^{(2)} : L^1 (\Gamma) \longrightarrow \Nh\Gamma \quad \mbox{with} \; \| r (f) \| \le \| f \|_1
\end{equation}
which extends the map (\ref{eq:4}) above. In particular, units in $L^1 (\Gamma)$ give units in $\Nh\Gamma$. For some classes of groups, the converse is known. We will discuss this in theorem \ref{t7} below.

{\bf Example} Take some element $g$ in $\C \Gamma$ with $\|g \|_1 < 1$. Then $1 + g$ is a unit in $L^1 (\Gamma)$ and hence $u = r (1+g) = 1 + r (g)$ is a unit in $\Nh\Gamma$. The series
\[
\log (1+g) = \sum^{\infty}_{\nu = 1} \frac{(-1)^{\nu-1}}{\nu} g^{\nu}
\]
converges in $L^1 (\Gamma)$. If $u$ is also positive e.g. if $g = h + h^* + hh^*$ for some $h \in \C\Gamma$ with $\|h \|_1 < \sqrt{2} - 1$ then we have by (\ref{eq:6})
\begin{equation}
  \label{eq:9}
  \log \ddet_{\Nh\Gamma} u = \tr_{\Nh\Gamma} \left( \sum^{\infty}_{\nu=1} \frac{(-1)^{\nu-1}}{\nu} g^{\nu} \right) = \sum^{\infty}_{\nu=1} \frac{(-1)^{\nu-1}}{\nu} \tr_{\Nh\Gamma} (g^{\nu}) \; .
\end{equation}
Recall that $g^{\nu}$ is the $\nu$-th power of $g$ in $\C\Gamma$ and $\tr_{\Nh\Gamma} (g^{\nu})$ its $e$-coefficient.

Next we look at the case $\Gamma = \Z^n$. The Fourier transform provides an isometric isomorphism of Hilbert spaces $\Fh : L^2 (\Z^n) \silo L^2 (T^n)$. On $L^1 \cap L^2$ the Fourier transform and its inverse are given by
\[
\Fh (c) (z) = \sum_{\nu \in \Z^n} c_{\nu} z^{\nu} \quad \mbox{and} \quad \Fh^{-1} (f)_{\nu} = \int_{T^n} f (z) \oz^{\nu} \, d\mu (z) \; .
\]
For an operator $A$ in $\Bh (L^2 (\Z^n))$ define the operator $\hA$ in $\Bh (L^2 (T^n))$ by the commutative diagram
\[
\xymatrix{
L^2 (\Z^n) \ar[r]^A \ar[d]^{\wr}_{\Fh} & L^2 (\Z^n) \ar[d]^{\wr \, \Fh} \\
L^2 (T^n) \ar[r]^{\hA} & L^2 (T^n) \; .
}
\]
Using the formulas for the Fourier transform if follows that for $\Z^n$-equivariant $A$ i.e. for $A$ in $\Nh \Z^n$ the operator $\hA$ is given by multiplication with the $L^{\infty}$-function $\Fh (A (0))$ on $T^n$. Here $0 \in \Z^n \subset L^2 (\Z^n)$. In this way one obtains an isomorphism of von Neumann algebras
\[
\Nh \Z^n \silo L^{\infty} (T^n) \; , \; A \longmapsto \Fh (A (0)) \; .
\]
On $L^{\infty} (T^n)$ the trace $\tr_{\Nh\Z^n}$ becomes the Haar integral $\int_{T^n}$ since
\[
\tr_{\Nh\Z^n} (A) = (A (0) , 0) = A (0)_0 = \int_{T^n} \Fh (A (0)) \, d\mu \; .
\]

\begin{prop}
  \label{t4}
For an element $f$ in $\C [\Z^n]$ the following conditions are equivalent:\\
{\bf a} \quad $f$ is a unit in $L^1 (\Z^n)$\\
{\bf b} \quad $f$ is a unit in $\Nh\Z^n$\\
{\bf c} \quad $f$ viewed as a Laurent polynomial does not vanish in any point of $T^n$.\\
Under these conditions we have
\[
\ddet_{\Nh\Gamma} f = \exp \int_{T^n} \log |f (z)| \, d\mu (z) \; .
\]
\end{prop}

\begin{proof}
  We have seen that {\bf a} implies {\bf b}. If $f$ is a unit in $\Nh\Z^n$ then $\Fh (r (f) (0))$ is a unit in $L^{\infty} (T^n)$. By definition, $r (f) (0) = f^* \in \C [\Z^n] \subset L^2 (\Z^n)$ and $\Fh (f^*) = \of$ where on the right $f$ is viewed as a Laurent polynomial in the coordinates $z = (z_1 , \ldots , z_n)$ on $T^n$. Since $\Fh (r (f) (0)) = \of$ is a continuous function on $T^n$ it is a unit in $L^{\infty} (T^n)$ if and only if it does not vanish in any point of $T^n$. Hence {\bf b} implies {\bf c}. The implication {\bf c} $\Rightarrow$ {\bf a} follows from a famous theorem of Wiener \cite{W} Lemma II\,e, who showed that if a continuous nowhere vanishing function $f$ on $T^n$ has Fourier coefficients in $L^1 (\Z^n)$, then $1 /f$ has Fourier coefficients in $L^1 (\Z^n)$ as well. See \cite{K} for a modern proof using Gelfand's theory of commutative $C^*$-algebras.

As we have seen, the element $f$ in $\C [\Z^n] \subset \Nh\Z^n$ maps to the function $\of$ in $L^{\infty} (T^n)$. Hence $\log ff^*$ in $\Nh\Z^n$ corresponds to $\log \of f = \log |f|^2$ in $L^{\infty} (T^n)$. Since the von Neumann trace on $\Nh \Z^n$ is the integral over $T^n$ the formula for $\det_{\Nh\Gamma} f$ follows.
\end{proof}

There is a classification theory for von Neumann algebras in terms of types. See \cite{L} 9.1.2 for a quick overview. For a finitely generated discrete group $\Gamma$ the von Neumann algebra $\Nh\Gamma$ is a type $\mathrm{I}_f$ if $\Gamma$ has an abelian subgroup of finite index. Otherwise it is of type $\mathrm{II}_1$. Thus for example $\Nh\Gamma$ is of type $\mathrm{II}_1$ for the $3$-dimensional integral Heisenberg group which is a first non-trivial example of a torsion-free non-commutative amenable discrete group.

We now recall a crucial approximation result due to L\"uck for von Neumann traces of amenable groups, \cite{L} Lemma 13.42. Let $\Gamma$ be a finitely generated discrete amenable group with a F{\o}lner sequence $(F_n)$. An element $A$ of $\Nh\Gamma$ induces for every finite subset $F \subset \Gamma$ an endomorphism, obtained by composition:
\begin{equation}
  \label{eq:10}
  A_F : \C [F] \overset{i_F}{\hookrightarrow} L^2 (\Gamma) \xrightarrow{A} L^2 (\Gamma) \overset{p_F}{\twoheadrightarrow} \C [F] \; .
\end{equation}
Here $i_F$ and $p_F$ are the canonical inclusion and projection maps respectively. Note that $p^*_F = i_F$ and hence $(A^*)_F = (A_F)^*$. For the operator norms adapted to the $L^2$-norms we have $\|A_F \| \le \|A \|$. We write $A_n = A_{F_n}$ etc.

\begin{theoremon}
  {\bf (L\"uck)} For every complex polynomial $Q$ and any $A$ in $\Nh\Gamma$ we have
\[
\tr_{\Nh\Gamma} (Q(A)) = \lim_{n\to\infty} \frac{1}{|F_n|} \tr (Q (A_n)) \; .
\]
\end{theoremon}

For the convenience of the reader let us sketch the {\it proof}. It suffices to take $Q (T) = T^s$ for $s \ge 1$. Let $\pi_n = i_n p_n$ be the orthogonal projection of $L^2 (\Gamma)$ to $\C [F_n]$. For $\gamma \in F_n \subset \C [F_n]$ we have:
\[
(A^s_n (\gamma), \gamma) = ((\pi_n A)^s (\gamma) , \gamma) \; .
\]
Using the telescope formula
\[
A^s - (\pi_n A)^s = \sum^{s-1}_{\nu=0} (\pi_n A)^{\nu} (1 - \pi_n) A^{s-\nu}
\]
one finds the estimate
\begin{equation}
  \label{eq:11}
  |(A^s (\gamma),\gamma) - (A^s_n (\gamma) , \gamma)| \le \sum^{s-1}_{\nu=0} \|(1 - \pi_n) A^{s-\nu} (\gamma)\|_2 \, \|A \|^{\nu} \; .
\end{equation}
For $\varepsilon > 0$, writing
\[
A^i (e) = \sum_{\gamma} \lambda^{(i)}_{\gamma} \gamma \quad \mbox{in} \; L^2 (\Gamma) \quad \mbox{for} \; 1 \le i \le s
\]
one finds an integer $R = R (\varepsilon) \ge 1$ such that we have:
\[
\Big\| \sum_{d (\gamma , e) > R} \lambda^{(i)}_{\gamma} \gamma \Big\|_2 < \varepsilon \quad \mbox{for all} \; 1 \le i \le s \; .
\]
Here $d$ is the left invariant metric on the (right) Cayley graph $\Gh$ of $\Gamma$ with respect to a finite set $S$ of generators of $\Gamma$ with $S = S^{-1}$. Its set of vertices is $\Gamma$ and two vertices $\gamma , \gamma'$ of $\Gh$ are connected by a (unique) edge if and only if $\gamma' = \gamma s$ for some $s$ in $S$. 

For all $n \ge 1$ we get:
\begin{eqnarray*}
  \Big| \tr_{\Nh\Gamma} (A^s) - \frac{1}{|F_n|} \tr (A^s_n) \Big| & = & \Big| \frac{1}{|F_n|} \sum_{\gamma \in F_n} (A^s (\gamma) , \gamma) - (A^s_n (\gamma) , \gamma) \Big| \\
& \le & \frac{1}{|F_n|} \sum_{\gamma \in F_n} \sum^{s-1}_{\nu=0} \| (1 - \pi_n) A^{s-\nu} (\gamma) \|_2 \, \|A\|^{\nu} \; .
\end{eqnarray*}
For $\gamma$ with $B_R (\gamma) \subset F_n$ we have:
\[
\| (1 - \pi_n) A^{s-\nu} (\gamma) \|_2 < \varepsilon
\]
since $\gamma\gamma' \notin F_n$ implies that $d (\gamma' , e) = d (\gamma \gamma' , \gamma) > R$. Separating the sum over $\gamma \in F_n$ into a sum over $\gamma$'s with $B_R (\gamma) \not\subset F_n$ and another one over the rest, one gets the estimate
\[
\Big| \tr_{\Nh\Gamma} (A^s) - \frac{1}{|F_n|} \tr (A^s_n) \Big| \le \frac{| \{ \gamma \in F_n \tei B_R (\gamma) \not\subset F_n \}|}{|F_n|} s \|A \|^s + s\varepsilon \max (1 , \|A \|^s) \; .
\]
Let $K$ be the finite set of $\le R$-fold products of element in $S$. Then $B_R (\gamma) = \gamma K$ and hence $B_R (\gamma) \not\subset F_n$ is equivalent to $\gamma \notin F_n K^{-1} = F_n K$. Thus by the F{\o}lner property the first term on the right tends to zero for $n \to \infty$. This gives the assertion. \endproof

In general we cannot approximate $\det_{\Nh\Gamma} (u)$ for $u$ in $(\Nh\Gamma)^*$ by the finite dimensional determinants $\det u_n$ because these may well be zero for all $n$. If for example $F_n \gamma^{-1} \not\subset F_n$, then for $u = r (\gamma)$ we have $\det_{\Nh\Gamma} (u) = 1$ and $\det u_n = 0$. If however $u$ is positive this cannot occur and we have the following result:

\begin{theorem}
  \label{t5}
If $u$ in $(\Nh\Gamma)^*$ is positive then all induced endomorphisms $u_n$ of $\C [F_n]$ are positive as well and we have
\[
\ddet_{\Nh\Gamma} u = \lim_{n\to \infty} (\ddet u_n)^{1 / |F_n|} \; .
\]
\end{theorem}

\begin{proof}
  For a finite subset $F$ of $\Gamma$ consider the endomorphism $u_F = p_F \, u \, i_F$ of $\C [F]$ as in formula (\ref{eq:10}). Since $u$ is selfadjoint, $u_F$ is selfadjoint as well. Moreover for $v \neq 0$ in $\C [F]$ we have since $u$ is positive:
  \begin{equation}
    \label{eq:12}
    (u_F (v) , v) = (p_F \, u \, i_F (v) , v) = (u  (i_F (v)) , i_F (v)) > 0 \; .
  \end{equation}
Hence $u_F$ is positive as well and in particular an isomorphism. Equation (\ref{eq:12}) implies that we have:
\[
\sup_{v \in \C [F] \atop \| v\|_2 = 1} (u_F (v) , v) \le \sup_{v \in L^2 (\Gamma) \atop \| v \|_2 = 1} (u (v) , v)
\]
and
\[
\inf_{v \in \C [F] \atop \| v\|_2 = 1} (u_F (v) , v) \ge \inf_{v \in L^2 (\Gamma) \atop \| v \|_2 = 1} (u (v) , v) \; .
\]
Now, according to \cite{Y} XI,8, Theorem 2 and VII,3, Theorem 3 we have the following information about the spectrum $\sigma (T)$ of a bounded selfadjoint operator $T$ on a Hilbert space:
\begin{equation}
\label{eq:13_neu}
\sup_{\lambda \in \sigma (T)} \lambda = \sup_{\| v \| = 1} (Tv , v) = \| T \| \quad \mbox{and} \quad \inf_{\lambda \in \sigma (T)} \lambda = \inf_{\| v \| = 1} (Tv , v) \; .
\end{equation}
Since $\sigma (u)$ is a compact subset of $(0,\infty)$ since $u$ is positive it follows that we have
\begin{equation}
  \label{eq:13}
  0 < a := \inf_{\lambda \in \sigma (u)} \lambda \le \inf_{\lambda \in \sigma (u_F)} \lambda \le \sup_{\lambda \in \sigma (u_F)} \lambda \le \sup_{\lambda \in \sigma (u)} \lambda =: b \; .
\end{equation}
The operators $\log u$ in $\Nh\Gamma$ and $\log u_n$ in $\End \C [F_n]$ can be defined by the functional calculus since $u$ and $u_n$ are positive. For the theorem, we have to prove the limit formula:
\[
\tr_{\Nh\Gamma} (\log u) = \lim_{n\to\infty} \frac{1}{|F_n|} \tr (\log u_n) \; .
\]
Fix $\varepsilon > 0$. By the Weierstra{\ss} approximation theorem there exists a real polynomial $Q$ such that
\[
\sup_{t \in [a,b]} |\log t - Q (t)| \le \varepsilon \; .
\]
By the estimates (\ref{eq:13}) we know that $\sigma (u), \sigma (u_n) \subset [a,b]$ for all $n \ge 1$. Hence we get the following estimates for the spectral norms i.e. the operator norms adapted to the $L^2$-norms:
\[
\| \log u - Q (u) \| \le \varepsilon \quad \mbox{and} \quad \|\log u_n - Q (u_n) \| \le \varepsilon \; .
\]
By L\"uck's approximation result above, we may choose $n_0 \ge 1$, such that for all $n \ge n_0$ we have
\[
\Big| \tr_{\Nh\Gamma} Q (u) - \frac{1}{|F_n|} \tr Q (u_n) \Big| \le \varepsilon \; .
\]
For any selfadjoint endomorphism $T$ of $\C [F]$ we have the estimate
\[
|\tr \, T | \le |F| \, \| T \|
\]
because
\[
\Big| \sum_{\lambda \in \sigma (T)} \lambda \Big| \le \dim \C [F] \sup_{\lambda \in \sigma (T)} \lambda = |F| \sup_{\| v \|_2 = 1} (Tv, v) \le |F|\, \| T \| \; .
\]
Hence we get:
\begin{eqnarray*}
  \lefteqn{\Big|\tr_{\Nh\Gamma} (\log u) - \frac{1}{|F_n|} \tr (\log u_n) \Big| \le |\tr_{\Nh\Gamma} (\log u - Q (u))| + } \\
& & + \Big| \tr_{\Nh\Gamma} Q (u) - \frac{1}{|F_n|} \tr Q (u_n)\Big| + \frac{1}{|F_n|} |\tr (\log u_n - Q (u_n))| \\[0.2cm]
& \le & \| \log u - Q (u)\| + \varepsilon + \| \log u_n - Q (u_n)\| \le 3 \varepsilon \; .
\end{eqnarray*}
\end{proof}

\begin{prop}
  \label{t7_neu}
For positive $u$ in $(\Nh\Gamma)^{\times}$ we have $\|u^{-1}_F \| \le \| u^{-1} \|$ for the spectral norms. 
\end{prop}

\begin{proof}
  This follows from equations (\ref{eq:13_neu}) and (\ref{eq:13}):
\[
\|u^{-1}_F\| = \sup_{\lambda \in \sigma (u_F)} \lambda^{-1} = \big( \inf_{\lambda \in \sigma (u_F)}\lambda \big)^{-1} \le \big( \inf_{\lambda \in \sigma (u)} \lambda \big)^{-1} = \sup_{\lambda \in \sigma (u)} \lambda^{-1} = \| u^{-1}\| \; .
\]
\end{proof}

For $f$ in $\C\Gamma$ and finite $F \subset \Gamma$ we denote the endomorphism $r (f)_F$ of $\C [F]$ also by $f_F$. It is equal to the composition
\[
f_F : \C [F] \overset{i}{\hookrightarrow} \C [\Gamma] \xrightarrow{R_{f^*}} \C [\Gamma] \overset{p}{\twoheadrightarrow} \C [F] \;.
\]
If $f$ is in $\Z \Gamma$, then $f_F$ maps $\Z [F]$ and $\R [F]$ into themselves. Since $\Z \Gamma$ is a lattice in $\R \Gamma$ we get the following corollary of theorem \ref{t5}, where we write $f_n$ for $f_{F_n}$:

\begin{cor}
  \label{t6}
Let $\Gamma$ be a finitely generated discrete amenable group and consider an element $f$ of $\Z \Gamma$ which is a positive unit in $\Nh\Gamma$. Then for every F{\o}lner sequence $(F_n)$ we have the formula:
\[
\ddet_{\Nh\Gamma} f = \lim_{n\to\infty} |\Z [F_n] / f_n \Z [F_n]|^{1/|F_n|} \; .
\]
\end{cor}

\begin{proof}
  For any isomorphism $\varphi$ of a finite dimensional real vector space $V$ which maps a lattice $\Lambda$ of $V$ into itself we have
\[
|\det \varphi| = |\Lambda / \varphi (\Lambda)| \; .
\]
Hence the corollary follows from theorem \ref{t5}.
\end{proof}

The next result follows by combining several results in the literature.

\begin{theorem}
  \label{t7}
Let $\Gamma$ be a finitely generated virtually nilpotent discrete group. Then we have: 

(W) \hfill \parbox{85ex}{\[
    L^1 (\Gamma) \cap (\Nh\Gamma)^{\times} = L^1 (\Gamma)^{\times} \; .
  \]}
\end{theorem}

\begin{proof}
  The one-sided Wiener property (W') for a discrete group $\Gamma$ is that every proper left (right) ideal in $L^1 (\Gamma)$ is annihilated by a non-zero positive linear functional, \cite{Le}. Equivalently, every such ideal has to annihilate a non-zero vector in a non-degenerate $*$-representation of $L^1 (\Gamma)$. Another way to phrase this condition is as follows:

{\bf Claim} Condition (W') is equivalent to the condition

(W'') \hfill \parbox{85ex}{\[
    L^1 (\Gamma) \cap C^* (\Gamma)^{\times} = L^1 (\Gamma)^{\times}
  \]}

where $C^* (\Gamma)$ is the $C^*$-algebra of $\Gamma$. 

\begin{proof}
Assume that $\Gamma$ satisfies (W') and consider $f$ in $L^1 (\Gamma) \cap C^* (\Gamma)^{\times}$. Then for every non-degenerate $*$-representation $\pi$ of $L^1 (\Gamma)$ or equivalently for every such representation of $C^* (\Gamma)$ the operator $\pi (f)$ is invertible. By condition (W') we must therefore have $L^1 (\Gamma) f = L^1 (\Gamma)$ and $fL^1 (\Gamma) = L^1 (\Gamma)$. Hence $f$ is a unit in $L^1 (\Gamma)$. Now assume that $\Gamma$ does not satisfy (W'). Then there is a proper left ideal $\ea$ in $L^1 (\Gamma)$ which is not annihilated by a state of $C^* (\Gamma)$. By \cite{Di} Theorem 2.9.5, it follows that $\ea$ must be dense in $C^* (\Gamma)$. Hence there is an element $f$ in $\ea \cap C^* (\Gamma)^{\times}$. It cannot be a unit in $L^1 (\Gamma)$ and hence (W'') is not satisfied. This proves the claim.
\end{proof}

Let $C^*_r(\Gamma)$ be the reduced $C^*$-algebra of $\Gamma$ i.e. the closure of $L^1 (\Gamma)$ in the operator norm of $\Bh (L^2 (\Gamma))$. Using the functional calculus one sees that
\[
C^*_r (\Gamma) \cap (\Nh\Gamma)^{\times} = C^*_r (\Gamma)^{\times} \; .
\]
If the discrete group $\Gamma$ is amenable, we have $C^*_r (\Gamma) = C^* (\Gamma)$ and hence (W) is equivalent to (W'') i.e. to (W'). It follows from the work of Losert \cite{Lo} that a discrete group satisfies (W') if and only if it is ``symmetric''. Finitely generated virtually nilpotent discrete groups are symmetric by \cite{LeP}, Corollary 3.
\end{proof}

\begin{rem}
  The argument shows that for finitely generated discrete amenable groups condition (W) is equivalent to the symmetry of $\Gamma$. 
\end{rem}
%\newpage
%\input{sec4_fuglede}

\section{The dynamical systems}
\label{sec:4}

In this section we introduce the class of dynamical systems that we are interested in and make a remark on expansiveness.

For a discrete group $\Gamma$ consider the Pontrjagin pairing
\[
\langle , \rangle : \Z \Gamma \times \widehat{\Z \Gamma} \longrightarrow \R / \Z \; .
\]
If we view the elements of $\widehat{\Z\Gamma} = \R / \Z [|\Gamma|]$ as infinite formal sums over $\Gamma$ with coefficients in $\R / \Z$, we have:
\[
\Big\langle \sum_{\gamma} a_{\gamma} \gamma , \sum_{\gamma} x_{\gamma} \gamma \Big\rangle = \sum_{\gamma} a_{\gamma} x_{\gamma} \quad \mbox{in} \; \R / \Z \; .
\]
Note that the sum on the right is finite. The ring $\Z \Gamma$ operates by left and right multiplication on the abelian group $\R / \Z [|\Gamma|]$. The following formula holds for all $a,f$ in $\Z \Gamma$ and $x$ in $\R / \Z [|\Gamma|]$:
\begin{equation}
  \label{eq:15}
  \langle af , x \rangle = \langle a , xf^* \rangle = \langle f , a^*x \rangle \; .
\end{equation}
This needs to be checked for $a,f,x$ in $\Gamma$ only, where it is clear. According to (\ref{eq:15}) the Pontrjagin dual to right multiplication by $f$ on $\Z \Gamma$ is right multiplication by $f^*$ on $\R / \Z [|\Gamma|]$. Hence we have:
\[
X_f := \widehat{\Z \Gamma / \Z \Gamma f} = \Ker (R_{f^*} : \R / \Z [|\Gamma|] \longrightarrow \R / \Z [|\Gamma|]) \; .
\]
We have a natural left $\Gamma$-action on $\widehat{\Z\Gamma}$ by sending the character $\chi$ to the character $\gamma \chi$ defined by $(\gamma \chi) (a) = \chi (\gamma^{-1} a)$. According to (\ref{eq:15}) this corresponds to left multiplication by $\gamma \in \Gamma$ on $\R / \Z [|\Gamma|]$. Since left and right multiplication commute this left $\Gamma$-action passes to $X_f$. It is clear that $\Gamma$ acts by (topological) automorphisms on the compact abelian group $X_f$.

Now assume that $\Gamma$ is a finitely generated discrete amenable group. Since $\Z\Gamma / \Z\Gamma f$ is countable the compact group $X_f$ is metrizable. In section \ref{sec:2} we have seen that the different notions of entropy for the $\Gamma$-action on $X_f$ coincide and we set:
\[
h_f := h_{\cover} = h_{\sep} = h_{\spa} = h_{\mu}
\]
where $\mu$ is the Haar measure on $X_f$. Let $\vartheta$ be the invariant metric on $\R / \Z$ induced from $\R$ i.e. 
\[
\vartheta (t \bmod \Z , t' \bmod \Z) = \min_{n\in \Z} |t - t' + n | \; .
\]
According to proposition \ref{t23} applied to $X_f$ and $G = \R / \Z$ we have:
\begin{equation}
  \label{eq:16}
  h_f = \lim_{\varepsilon \to 0} \varlimsup_{n\to\infty} \frac{1}{|F_n|} \log s'_{F_n} (\varepsilon) \; .
\end{equation}
Here $s'_F (\varepsilon)$ is the maximal cardinality of a $[F,\varepsilon]$-separated subset $E$ of $X_f$. Recall, that $E$ is $[F,\varepsilon]$-separated if for all $x \neq y$ in $E$ there is some $\gamma$ in $F$ with $\vartheta (x_{\gamma} , y_{\gamma}) \ge \varepsilon$. This is the formula for the entropy $h_f$ that will be used later. 

Recall that a $\Gamma$-action by homeomorphisms on a compact metrizable space $X$ is expansive if there exists a metric $d$ defining the topology and an $\varepsilon > 0$ such that for all $x \neq y$ in $X$ we have $d (\gamma x , \gamma y) \ge \varepsilon$ for some $\gamma$ in $\Gamma$. If $\Gamma$ acts by automorphisms on a compact metrizable abelian topological group $X$, this is equivalent to the existence of an expansive neighborhood $N$ of $0 \in X$, i.e. $\bigcap_{\gamma\in\Gamma} \gamma N = 0$.

\begin{prop}
  \label{t9}
Let $\Gamma$ be a finitely generated discrete group. If $f$ in $\Z \Gamma$ is a unit in $L^1 (\Gamma)$ the $\Gamma$-action on $X_f$ is expansive.
\end{prop}

\begin{proof}
  Set $\varepsilon = 1/3 \| f \|_1$ and $N = \{ x \in X_f \tei \vartheta (x_e , 0) < \varepsilon \}$. We claim that $N$ is an expansive neighborhood of $0 \in X_f$. If not, there exists some $x \neq 0$ in $\bigcap_{\gamma \in \Gamma} \gamma N$ i.e. with $\vartheta (x_{\gamma} , 0) < \varepsilon$ for all $\gamma \in \Gamma$. Because of $\varepsilon \le 1/3$ there exists a unique $y$ in $L^{\infty} (\Gamma)$ with $|y_{\gamma}| < \varepsilon$ and $x_{\gamma} = y_{\gamma} \bmod \Z$ for all $\gamma \in \Gamma$. By definition of $X_f$ we have $R_{f^*} (x) = 0$. Hence $r^{(\infty)} (f) (y) = yf^* \equiv 0 \bmod \Z$ and therefore $\| r^{(\infty)} (f) (y)\|_{\infty} \ge 1$ if $r^{(\infty)} (f) (y)$ is nonzero. Estimate (\ref{eq:7}) on the other hand gives:
\[
\| r^{(\infty)} (f) (y) \|_{\infty} \le \| r^{(\infty)} (f) \| \, \|y \|_{\infty} \le \| f \|_1 \, \|y \|_{\infty} \le \| f \|_1 \cdot \varepsilon = 1/3
\]
and therefore $r^{(\infty)} (f) y = 0$. By assumption $f$ is a unit in $L^1 (\Gamma)$. Hence $r^{(\infty)} (f)$ is an automorphism of $L^{\infty} (\Gamma)$ and $y$ must be zero. This contradiction proves the proposition.
\end{proof}

\begin{rems}
  The proof mimics the first half of the proof of Lemma 6.8 in \cite{Sch}. For finitely generated virtually nilpotent groups it follows from theorem \ref{t7} and the proposition that $X_f$ is expansive if $f$ is a unit in $\Nh\Gamma$.
\end{rems}
%\newpage
%\input{sec5_fuglede}

\section{Proof of the inequality $h_f \ge \log \det_{\Nh\Gamma} f$}
\label{sec:5}

In this section we prove the following result:

\begin{theorem}
  \label{t10}
Let $\Gamma$ be a finitely generated discrete amenable group. If $f$ in $\Z\Gamma$ is a unit in $L^1 (\Gamma)$ and positive in $\Nh\Gamma$ then we have:
\[
h_f \ge \log \ddet_{\Nh\Gamma} f \; .
\]
\end{theorem}

\begin{proof}
  For a finite subset $F \subset \Gamma$ consider the isomorphism c.f. \S\,\ref{sec:3}
\[
f_F : \R [F] \xrightarrow{i} \R \Gamma \xrightarrow{R_{f^*}} \R \Gamma \xrightarrow{p} \R [F] \; .
\]
It maps $\Z [F]$ into itself. Setting
\[
Y_F = f^{-1}_F \Z [F] / \Z [F] \subset \R [F] / \Z [F] = \R / \Z [F]
\]
we have by corollary \ref{t6}
\begin{equation}
  \label{eq:18}
  \log \ddet_{\Nh\Gamma} f = \lim_{n\to\infty} \frac{1}{|F_n|} \log |Y_{F_n}| \; .
\end{equation}
For $t$ in $\R / \Z$ let $\tilde{t}$ in $\R$ be its unique representative with $0 \le \tilde{t} < 1$.

Let $\R / \Z [F] \to \R [F] , y \mapsto \tilde{y}$ be the map defined by $\ty_{\gamma} = (\widetilde{y_{\gamma}})$ for all $\gamma \in F$. For an element $z$ of $\R [|\Gamma |]$ let $[z]$ be its class in $\R / \Z [|\Gamma |]$. Since $f$ is invertible in $\Nh\Gamma$ the operator
\[
r (f)^{-1} : L^2 (\Gamma) \longrightarrow L^2 (\Gamma)
\]
exists and is bounded. In fact by (\ref{eq:8}) we have $\|r (f)^{-1} \| \le \| f^{-1} \|_1$ but this is not essential. The operator $r (f)^{-1}$ restricts to the real $L^2$-spaces and we can consider the composition:
\[
\R [F] \xrightarrow{f_F} \R [F] \overset{i_F}{\hookrightarrow} L^2 (\Gamma)_{\R} \xrightarrow{r (f)^{-1}} L^2 (\Gamma)_{\R} \hookrightarrow \R [|\Gamma |] \; .
\]

{\bf Claim} Setting $\varphi (y) = [r (f)^{-1} f_F (\ty)]$ we obtain a map $\varphi : Y_F \longrightarrow X_f$.

Namely, we have:
\[
R_{f^*} [r (f)^{-1} f_F (\ty)] = [r (f) r (f)^{-1} f_F (\ty)] = [f_F (\ty)] = 0
\]
since $f_F (\ty) \in \Z [F] \subset \Z [|\Gamma |]$ by the definition of $Y_F$. 

Now, the idea is to show that the map $\varphi$ injects $Y_F$ into an $[F', \varepsilon]$-separated subset of $X_f$ for suitable $\varepsilon > 0$ and a slight enlargement $F'$ of $F$. This does not work directly since $\varphi$ is not injective in general. However, using that $f$ is a unit in $L^1 (\Gamma)$ and not only in $\Nh\Gamma$ it is possible to replace $Y_F$ by a subset $\tY_F$ for which the above idea works. 

Let $K$ be the support of $f$ and set $A = \| f^{-1} \|_1 \| f \|_1$. Subdivide the interval $[-A , A]$ into disjoint subintervals $I_j , j \in J$ such that all $I_j$ with at most one exception have length $1/3$ and the remaining one has length $\le 1/3$. The inequality
\[
2A = \sum_{j\in J} |I_j| \ge \frac{1}{3} (|J| - 1)
\]
shows that $|J| \le 6A + 1$.

We claim that for all $y$ in $Y_F$ the number $r (f)^{-1} f_F (\ty)_{\gamma}$ lies in $[-A, A]$. Namely, since $f$ is invertible in $L^1 (\Gamma)$ the operator $r^{(\infty)} (f)$ is invertible and we can consider its inverse
\[
r^{(\infty)} (f)^{-1} : L^{\infty} (\Gamma) \longrightarrow L^{\infty} (\Gamma) \; .
\]
According to (\ref{eq:7}) we have $\| r^{\infty} (f)^{-1} \| \le \| f^{-1} \|_1$. Since $\R [F]$ lies in $L^2 (\Gamma) \subset L^{\infty} (\Gamma)$ we get:
\begin{eqnarray*}
  \| r (f)^{-1} f_F (\ty) \|_{\infty} & = & \| r^{(\infty)} (f)^{-1} f_F (\ty) \|_{\infty} \le \| r^{(\infty)} (f)^{-1}\| \, \| f_F (\ty) \|_{\infty} \\
& \le & \| f^{-1} \|_1 \, \| f \|_1 \, \| \ty \|_{\infty} \le A \; .
\end{eqnarray*}
Note here that the projection $L^{\infty} (\Gamma) \to L^{\infty} (F)$ has norm $\le 1$ and $\|r^{(\infty)} (f) \| \le \| f \|_1$.

I learned about the following type of argument from Rita Solomyak's paper \cite{S} \S\,3. We define an equivalence relation $\sim$ on $Y_F$ as follows. For elements $y , y'$ in $Y_F$ set $y \sim y'$ if for every $\gamma \in FK \ohne F$ there is some $j \in J$ such that both
\[
r (f)^{-1} f_F (\ty)_{\gamma} \quad \mbox{and} \quad r (f)^{-1} f_F (\ty')_{\gamma} \quad \mbox{lie in} \; I_j \; .
\]
(The above argument was needed to show that $\sim$ is reflexive.)

There are at most $|J|^{|FK \ohne F|} \le (6A + 1)^{|FK \ohne F|}$ equivalence classes. Fix an equivalence class $\tY_F$ of maximal order. Then we have the estimates
\[
|\tY_F | \le |Y_F| \le (6A + 1)^{|FK \ohne F|} |\tY_F| \; .
\]
Applying this to $F = F_n$ for a F{\o}lner sequence $(F_n)$ we get
\[
\frac{1}{|F_n|} \log |\tY_{F_n}| \le \frac{1}{|F_n|} \log |Y_{F_n}| \le \frac{1}{|F_n|} \log |\tY_{F_n}| + \frac{|F_n K \ohne F_n|}{|F_n|} \log (6A + 1) \; .
\]
Using the F{\o}lner property and equation (\ref{eq:18}) this gives:
\[
\varlimsup_{n\to\infty} \frac{1}{|F_n|} \log |\tY_{F_n}| \le \log \ddet_{\Nh\Gamma} f \quad \mbox{and} \quad \log \ddet_{\Gamma} f \le \varliminf_{n\to\infty} \frac{1}{|F_n|} \log |\tY_{F_n}| \; .
\]
Hence we have
\begin{equation}
  \label{eq:19}
  \log \ddet_{\Nh\Gamma} f = \lim_{n\to\infty} \frac{1}{|F_n|} \log |\tY_{F_n}| \; .
\end{equation}

{\bf Claim} The map $\varphi : \tY_F \to X_f$ is injective and maps $\tY_F$ into an $[FK , \varepsilon]$-separated set whenever $\varepsilon < 1 / \| f \|_1$.

We have to show that for $y , y'$ in $\tY_F$ the inequalities
\begin{equation}
  \label{eq:20}
  \vartheta ([r (f)^{-1} f_F (\ty)]_{\gamma} \; , \; [r (f)^{-1} f_F (\ty')]_{\gamma}) < \varepsilon \quad \mbox{for} \; \gamma \; \mbox{in} \; FK
\end{equation}
imply that $y = y'$. They are equivalent to the equations
\begin{equation}
  \label{eq:21}
  r (f)^{-1} f_F (\ty - \ty')_{\gamma} = \alpha_{\gamma} + \mu_{\gamma} \quad \mbox{for all} \; \gamma \; \mbox{in} \; FK
\end{equation}
with $\alpha_{\gamma} \in \Z$ and $|\mu_{\gamma}| < \varepsilon$.

By definition of $\tY_F$ we know that for every $\gamma$ in $FK \ohne F$ the real numbers $r (f)^{-1} f_F (\ty)_{\gamma}$ and $r (f)^{-1} f_F (\ty')_{\gamma}$ lie in the same interval $I_j$. Since $|I_j| \le 1/3$ it follows that we have
\[
|r (f)^{-1} f_F (\ty - \ty')_{\gamma}| \le 1/3 \; .
\]
Since $\varepsilon < 1/3$ we get for $\gamma$ in $FK \ohne F$:
\[
1/3 \ge |r (f)^{-1} f_F (\ty - \ty')_{\gamma}| = |\alpha_{\gamma} + \mu_{\gamma}| \ge |\alpha_{\gamma}| - |\mu_{\gamma}| \ge |\alpha_{\gamma}| - \varepsilon \ge |\alpha_{\gamma}| - 1/3 \; .
\]
Hence we have $|\alpha_{\gamma}| \le 2/3$ and therefore $\alpha_{\gamma} = 0$ since $\alpha_{\gamma} \in \Z$.

We can therefore write:
\begin{equation}
  \label{eq:22}
  r (f)^{-1} f_F (\ty - \ty') = \mu + \alpha
\end{equation}
where $\mu \in L^2 (\Gamma)_{\R}$ with $|\mu_{\gamma}| < \varepsilon$ for all $\gamma \in FK$ and $\alpha \in \Z [F]$. Applying $r (f)$ and the projection $p_F : L^2 (\Gamma)_{\R} \to \R [F]$ we get
\begin{equation}
  \label{eq:23}
  f_F (\ty - \ty') = p_F r (f) \mu + p_F r (f) \alpha = p_F r (f) \mu + f_F (\alpha) \; .
\end{equation}
Now note that
\begin{equation}
  \label{eq:24}
  p_F r (f) \mu = p_F r (f) p_{FK} \mu \; .
\end{equation}
Namely, writing $f = \sum_{\gamma'' \in K} a_{\gamma''} \gamma''$ and hence $f^* = \sum_{\gamma'' \in K} \alpha_{\gamma''} \gamma^{''-1}$ we have for all $\gamma$ in $F$ that
\[
(r (f) \mu)_{\gamma} = (\mu f^*)_{\gamma} = \sum_{\gamma' \gamma^{''-1} = \gamma \atop \gamma'' \in K} \mu_{\gamma'} a_{\gamma''} \; .
\]
This depends only on the components $\mu_{\gamma'}$ of $\mu$ with $\gamma' \in FK$. 

Using (\ref{eq:24}) we get:
\begin{equation}
  \label{eq:25}
  f_F (\ty - \ty' - \alpha) = p_F r (f) \mu'
\end{equation}
with $\mu' = p_{FK} \mu \in I_{\varepsilon} [FK]$, where $I_{\varepsilon} = (-\varepsilon , \varepsilon)$. Next we observe that
\[
\| p_F r (f) \mu'\|_{\infty} \le \| r^{(\infty)} (f) \mu' \|_{\infty} \le \| f \|_1 \, \| \mu' \|_{\infty} \; .
\]
This gives
\[
p_F r (f) \mu'\in I_{\varepsilon \| f \|_1} [F] \subset I_1 [F] \; .
\]
On the other hand, by the definition of $Y_F$ we know that $f_F (\ty - \ty' - \alpha)$ is an element of $\Z [F]$.

Since $I_1 [F] \cap \Z [F] = 0$, equation (\ref{eq:25}) shows that $f_F (\ty - \ty' - \alpha) = 0$ and hence that $\ty = \ty' + \alpha$ because $f_F$ is an automorphism of $\R [F]$. Thus we get $y = [\ty] = [\ty' + \alpha] = y'$ and the claim is proved.

Next, using the notations from formula (\ref{eq:16}), the claim proves that for $\varepsilon < 1 / \| f \|_1$ we have
\begin{equation}
  \label{eq:26}
  s'_{FK} (\varepsilon) \ge |\tY_F| \; .
\end{equation}
Using formula (\ref{eq:19}) this shows that
\[
\varlimsup_{n\to\infty} \frac{1}{|F_n|} \log s'_{F_n K} (\varepsilon) \ge \lim_{n\to \infty} \frac{1}{|F_n|} \log |\tY_{F_n}| = \log \ddet_{\Nh\Gamma} f \; .
\]
Easy estimates using the inequality
\[
|A \bigtriangleup C| \le |A \bigtriangleup B| + |B \bigtriangleup C|
\]
for finite subsets $A, B, C$ of a given set show that $(F_n K)$ is a F{\o}lner sequence and that $\lim_{n\to\infty} |F_n K| / |F_n| = 1$. It follows that for $\varepsilon < 1 / \| f \|_1$ we have
\[
\varlimsup_{n\to\infty} \frac{1}{|F_n K|} \log s'_{F_n K} (\varepsilon) \ge \log \ddet_{\Nh\Gamma} f
\]
and hence using formula (\ref{eq:16}) we get the desired estimate
\[
h_f \ge \log \ddet_{\Nh\Gamma} f \; .
\]
\end{proof}
%\newpage
%\input{sec6_fuglede}

\section{Proof of the inequality $h_f \le \log \det_{\Nh\Gamma} f$}
\label{sec:6}

Let $\Gamma$ be a discrete group. A sequence $(F_n)$ of finite subsets of $\Gamma$ will be called a ``strong F{\o}lner sequence'' if for every finite subset $K$ of $\Gamma$ we have
\begin{equation}
  \label{eq:27}
  \lim_{n\to\infty} \frac{|F_n K \ohne F_n|}{|F_n|} \log (1 + |F_n K \ohne F_n|) = 0 \; .
\end{equation}
This is implied by the following condition for all $\gamma$ in $\Gamma$:
\begin{equation}
  \label{eq:28}
  \lim_{n\to\infty} \frac{|F_n \gamma \ohne F_n|}{|F_n|} \log |F_n| = 0 \; .
\end{equation}
Note that it suffices to verify (\ref{eq:28}) for $\gamma$'s running over a symmetric set of generators of $\Gamma$.

\begin{prop}
  \label{t11}
Let $\Gamma$ be finitely generated and let $S$ be a finite system of generators with $e \in S = S^{-1}$.\\
1) If the condition
\begin{equation}
  \label{eq:29}
  \varliminf_{n\to\infty} \left( \frac{|S^{n+1}|}{|S^n|} - 1 \right) \log |S^n| = 0
\end{equation}
holds, then a suitable subsequence of $(S^n)_{n \ge 1}$ is a strong F{\o}lner sequence.\\
2) Condition (\ref{eq:29}) is satisfied if for every $\varepsilon > 0$ we have
\begin{equation}
  \label{eq:30}
  |S^n| \le \exp (\varepsilon \sqrt{n}) \quad \mbox{for all} \; n \ge n (\varepsilon) \; .
\end{equation}
In particular, virtually nilpotent groups being of polynomial growth have strong F{\o}lner sequences.
\end{prop}

\begin{rems}
Working with (\ref{eq:27}) instead of (\ref{eq:28}) does not allow faster growth rates than (\ref{eq:30}).  It does not seem to be known whether groups not of polynomial growth can have the growth rate (\ref{eq:30}). See \cite{Gri} for a survey on growth in group theory. On the other hand, it is not clear to me whether there are finitely generated amenable groups which do not have any strong F{\o}lner sequence.
\end{rems}

\begin{proof}
  1) Condition (\ref{eq:29}) implies that a subsequence of $(S^n)_{n \ge 1}$ satisfies (\ref{eq:28}) for all $\gamma$ in $S$.\\
2) If (\ref{eq:29}) does not hold, then setting $a_n = |S^n|$, there is some $\delta > 0$ such that for all $n \gg 0$ we have
\[
\left( \frac{a_{n+1}}{a_n} - 1 \right) \log a_n \ge \delta \; .
\]
This implies that for some constant $c > 0$ we have for $n \gg 0$
\[
\log a_{n+1} - \log a_n \ge \log \left( 1 + \frac{\delta}{\log a_n} \right) \ge \frac{c}{\log a_n} \; .
\]
Setting $b_n = \log a_n$ we find
\[
b^2_{n+1} - b^2_n \ge b_{n+1} b_n - b^2_n \ge c
\]
and hence $b_n \ge \lambda \sqrt{n}$ for some $\lambda > 0$ and all $n \gg 0$. This contradicts condition (\ref{eq:30}).
\end{proof}

\begin{theorem}
  \label{t12}
Let $\Gamma$ be a finitely generated discrete group which has a strong F{\o}lner sequence and hence is amenable. If $f$ in $\Z\Gamma$ is a positive unit in $\Nh\Gamma$ then we have
\[
h_f \le \log \ddet_{\Nh\Gamma} f \; .
\]
\end{theorem}

\begin{rem}
  For finite $F \subset \Gamma$ consider the isomorphism $f_F : L^{\infty} (F) \to L^{\infty} (F)$. The proof below shows that we do not have to assume the existence of a {\it strong} F{\o}lner sequence if the operator norms $\| f^{-1}_F \|$ adapted to the $L^{\infty}$-norms are uniformly bounded in $F$. As it is, we can only show this for the operator norms adapted to the $L^2$-norms. To make up for this weaker information the strong F{\o}lner condition is needed.
\end{rem}

\begin{proof}
  Consider the composition 
  \begin{equation}
    \label{eq:31}
    \R [| \Gamma |] \xrightarrow{R_{f^*}} \R [| \Gamma |] \xrightarrow{p_F} \R [F] \xrightarrow{f^{-1}_F} \R [F] \; .
  \end{equation}
It is a left inverse to the map $r (f)^{-1} i_F f_F$ used in the last section.

Let $\R / \Z [|\Gamma|] \to \R [|\Gamma|] , x \mapsto \tx$ be the map defined by $\tx_{\gamma} = (\widetilde{x_{\gamma}})$ for all $\gamma \in \Gamma$. 

{\bf Claim} Setting $\psi (x) = [f^{-1}_F p_F R_{f^*} (\tx)]$ we obtain a map $\psi : X_f \to Y_F$. 

Namely, since $R_{f^*} (x) = 0$ in $\R / \Z [|\Gamma|]$ we have:
\[
f_F (f^{-1}_F p_F R_{f^*} (\tx)) = p_F R_{f^*} (\tx) \in p_F \Z [|\Gamma|] = \Z [F] \; .
\]

Now, the idea is to show that for every $\varepsilon > 0$ and all $n \ge n_0 (\varepsilon)$ the map $\psi$ injects a maximal $[F_n , \varepsilon]$-separated subset $E$ of $X_f$ into $Y_{F_n}$. Again this does not work directly but only after replacing $E$ by a maximal equivalence class $\tilde{E}$ with respect to a suitable equivalence relation c.f. \cite{S} \S\,3.

Let $E [F,\varepsilon]$ be a maximal $[F,\varepsilon]$-separated subset of $X_f$. For any $L > 0$ we define an equivalence relation on $E [F,\varepsilon]$ as follows. Subdivide the interval $[0,1)$ into disjoint subintervals $I_j , j \in J$ such that all $I_j$ with at most one exception have length $L^{-1}$ and the remaining one has length $\le L^{-1}$. Then we have $|J| \le 1 + L$. 

We set $x \sim x'$ for elements $x, x'$ in $E [F,\varepsilon]$, if for every $\gamma \in FK \ohne F$ there is some $j \in J$ such that both $\tx_{\gamma}$ and $\tx'_{\gamma}$ lie in $I_j$. Here $K$ is a the support of $f$. There are at most $|J|^{|FK \ohne F|} \le (1+L)^{|FK\ohne F|}$ equivalence classes. Fix an equivalence class $\tilde{E} [F,\varepsilon]$ of maximal order. Then we have the estimates
\begin{equation}
  \label{eq:32}
  |\tilde{E} [F,\varepsilon]| \le |E [F, \varepsilon]| \le (1+L)^{|FK\ohne F|} |\tilde{E} [F,\varepsilon]| \; .
\end{equation}
In the following, the choice of $L$ will depend on $\varepsilon$ and $F$ i.e. $L = L_F (\varepsilon)$. For a F{\o}lner sequence $(F_n)$ set $L_n (\varepsilon) = L_{F_n} (\varepsilon)$. Then we have:
\begin{equation}
  \label{eq:33}
  \frac{1}{|F_n|} \log |\tilde{E} [F_n, \varepsilon]| \le \frac{1}{|F_n|} \log | E [F_n , \varepsilon]| \le \frac{1}{|F_n|} \log |\tilde{E} [F_n , \varepsilon]| + \frac{|F_nK \ohne F_n|}{|F_n|} \log (1 + L_n (\varepsilon)) \; .
\end{equation}
Now assume the condition
\begin{equation}
  \label{eq:34}
  \lim_{n\to\infty} \frac{|F_n K \ohne F_n|}{|F_n|} \log (1 + L_n (\varepsilon)) = 0 \quad \mbox{for every} \; \varepsilon > 0 \; .
\end{equation}
Then it follows from (\ref{eq:33}) and (\ref{eq:16}) that we have
\begin{equation}
  \label{eq:35}
  h_f = \lim_{\varepsilon \to 0} \varlimsup_{n \to\infty} \frac{1}{|F_n|} \log |\tilde{E} [F_n , \varepsilon]| \; .
\end{equation}
\end{proof}

{\bf Claim} For $\varepsilon > 0$, setting $\displaystyle L_F (\varepsilon) = \frac{2}{\varepsilon} \| r(f)^{-1} \| \, \|r (f) \| \, |FK \ohne F|^{1/2}$, the map \\
$\psi : \tilde{E} [F , \varepsilon] \to Y_F$ is injective.

\begin{proof}
  As in (\ref{eq:24}) the following equations hold for every $\mu$ in $\R [|\Gamma|]$:
  \begin{equation}
    \label{eq:36}
    p_F R_{f^*} (\mu) = p_F R_{f^*} p_{FK} (\mu) = p_F r(f) p_{FK} (\mu) \; .
  \end{equation}
For elements $x,x'$ in $\tilde{E} [F,\varepsilon]$ with $\psi (x) = \psi (x')$ we therefore have:
\[
f^{-1}_F p_F r (f) p_{FK} (\tx - \tx') \equiv 0 \bmod \Z [F] \; .
\]
Writing $p_{FK} (\tx - \tx') = p_F (\tx - \tx') + p_{FK \ohne F} (\tx - \tx')$ we get:
\[
p_F (\tx) - p_F (\tx') \equiv -f^{-1}_F p_F r (f) p_{FK \ohne F} (\tx - \tx') \bmod \Z [F] \; .
\]
Hence
\[
\max_{\gamma \in F} \vartheta (x _{\gamma} , x'_{\gamma}) \le \| f^{-1}_F p_f r(f) p_{FK \ohne F} (\tx - \tx') \|_{\infty} \; .
\]
It is at this point that we need an estimate for $\|f^{-1}_F z \|_{\infty}$. If there is some constant $c> 0$ such that $\| f^{-1}_F z \|_{\infty} \le c \|z \|_{\infty}$ for all $F$ and all $z$ in $L^{\infty} (F)$, then the constant $L = L_F (\varepsilon)$ can be chosen to be independent of $F$ and (\ref{eq:34}) is valid for any F{\o}lner sequence. Incidentially, at the corresponding stage of the argument in \cite{S} \S\,3 Rita Solomyak uses the maximum principle for a graph Laplacian. This gives good $L^{\infty}$-estimates directly. However, since we do not have such $L^{\infty}$-estimates, we have to be content with the estimate for $L^2$-norms of proposition \ref{t7_neu}:
\[
\|f^{-1}_F z \|_2 \le \| r (f)^{-1}\| \, \|z \|_2 \; .
\]
We get:
\begin{eqnarray*}
\max_{\gamma \in F} \vartheta (x_{\gamma} , x'_{\gamma}) & \le & \| f^{-1}_F p_F r (f) p_{FK \ohne F} (\tx - \tx') \|_2 \\
 & \le & \| r (f)^{-1} \| \, \|r (f)\| \, \| p_{FK \ohne F} (\tx - \tx')\|_2 \; .
\end{eqnarray*}
Now, for $w$ in $\R [FK \ohne F]$ we have the estimate:
\[ 
\| w \|_2 \le |FK \ohne F|^{1/2} \| w \|_{\infty} \; .
\]
This implies that
\[
\max_{\gamma \in F} \vartheta (x_{\gamma} , x'_{\gamma}) \le \frac{\varepsilon}{2} L_F (\varepsilon) \| p_{FK \ohne F} (\tx) - p_{FK \ohne F} (\tx') \|_{\infty} \; .
\]
By definition of the equivalence relation using $L = L_F (\varepsilon)$ on $E [F, \varepsilon]$ we know that since $x$ and $x'$ are equivalent:
\[
\| p_{FK \ohne F} (\tx) - p_{FK \ohne F} (\tx') \|_{\infty} \le L_F (\varepsilon)^{-1} \; .
\]
It follows that we have
\[
\max_{\gamma \in F} \vartheta (x_{\gamma} , x'_{\gamma}) \le \frac{\varepsilon}{2} < \varepsilon \; .
\]
Hence $x = x'$ since the set $\tilde{E} [F, \varepsilon] \subset E [F, \varepsilon]$ is $[F, \varepsilon]$-separated. Thus $\psi$ is injective on $\tilde{E} [F, \varepsilon]$ and the claim is proved.

Now the theorem is proved as follows. Let $(F_n)$ be a strong F{\o}lner sequence. Then with $L_F (\varepsilon)$ as in the claim, condition (\ref{eq:34}) is satisfied. Hence the entropy formula (\ref{eq:35}) holds with sets $\tilde{E} [F_n, \varepsilon]$ that inject into $Y_{F_n}$, according to the claim.

Using corollary \ref{t6}, we find:
\[
h_f \le \lim_{\varepsilon \to 0} \varlimsup_{n\to\infty} \frac{1}{|F_n|} \log |Y_{F_n}| = \log \ddet_{\Nh\Gamma} f \; .
\]
\end{proof}

\section{Further remarks}
\label{sec:7}

The first remark concerns non-expansive systems $X_f$. For L\"uck's generalization of the Fuglede--Kadison determinant, equation (\ref{eq:1}) is in fact valid for all $f \neq 0$ in $\Z [\Z^n]$. By the result of Lind, Schmidt and Ward \cite{LSW} the equation
\begin{equation}
  \label{eq:37}
  h_f = \log \ddet_{\Nh\Gamma} f
\end{equation}
therefore holds for all nonzero $f$ if $\Gamma = \Z^n$. However, for finite groups $\Gamma$ it does not hold for a single $f$ in $\Z \Gamma$ which is not a unit in $\Nh\Gamma$. Namely, according to \cite{L}, Example 3.12 we have
\begin{equation}
  \label{eq:38}
  \log \ddet_{\Nh\Gamma} f = \frac{1}{2|\Gamma|} \sum_{\lambda} \log \lambda \; .
\end{equation}
Here $\lambda$ runs over the nonzero eigenvalues of the selfadjoint positive semidefinit endomorphism $R_{ff^*}$ of $\C \Gamma$. On the other hand we have
\begin{equation}
  \label{eq:39}
  h_f = \sup_{\eU} \frac{1}{|\Gamma|} \log N (\eU^{\Gamma}) \; ,
\end{equation}
where $\eU$ runs over the open coverings of $X_f$. The space $X_f = \widehat{\Z\Gamma / \Z\Gamma f}$ consists of finitely many points if and only if $R_f (\Z \Gamma)$ has finite index in $\Z \Gamma$. This is equivalent to $R_f$ being an isomorphism of $\C\Gamma$ i.e. to $f^*$ and hence $f$ being a unit in $\Nh\Gamma$. If $f$ is not a unit in $\Nh\Gamma$ then the connected component of zero of $X_f$ is a real torus of positive dimension. From formula (\ref{eq:39}) we therefore get:
\begin{equation}
  \label{eq:40}
  h_f = \frac{1}{|\Gamma|} \log |X_f| = \frac{1}{|\Gamma|} \log |\Z \Gamma / \Z \Gamma f|
\end{equation}
if $f \in (\Nh\Gamma)^{\times}$ and $h_f = \infty$ if $f \notin (\Nh\Gamma)^{\times}$. Thus, for $f \in (\Nh\Gamma)^{\times}$ we have
\[
h_f = \log \ddet_{\Nh\Gamma} f
\]
in accordance with theorem \ref{t1_neu}. For $f \notin (\Nh\Gamma)^{\times}$ on the other hand we have the strict inequality
\[
h_f = \infty > \log \ddet_{\Nh\Gamma} f \; .
\]
Of course, the case of finite groups $\Gamma$ is somewhat exceptional for entropy theory and it could be that e.g. for torsion free groups $\Gamma$ equation (\ref{eq:38}) holds for all $f \neq 0$ in $\Z \Gamma$. 

We now comment on condition c) in theorem \ref{t1_neu}. 
The Yuzvinskii addition formula asserts (if true) that entropy is additive in short exact sequences of compact metrizable abelian groups with $\Gamma$-action. For $\Gamma = \Z^n$ this is known, \cite{Sch}, Theorem 14.1.

Let $f,g$ be elements of $\Z \Gamma \cap (\Nh\Gamma)^{\times}$. Then one has a canonical exact sequence of discrete abelian groups with $\Gamma$-action:
\[
0 \longrightarrow \Z \Gamma / \Z \Gamma f \xrightarrow{R_g} \Z \Gamma / \Z \Gamma fg \longrightarrow \Z \Gamma / \Z \Gamma g \longrightarrow 0 \, .
\]
Note that $R_g$ is injective because $g$ is not a zero divisor in $\Nh \Gamma \supset \Z \Gamma$. Dualizing and applying a Yuzvinskii addition formula would give $ h_{fg} = h_f + h_g$.
For $f$ in $\Z \Gamma \cap (\Nh\Gamma)^{\times}$ the element $f^*$ lies in $\Z \Gamma \cap (\Nh\Gamma)^{\times}$ as well and we expect that $h_f = h_{f^*}$. This seems to be a nontrivial relation because in general $X_f$ and $X_{f^*}$ are not isomorphic as topological {\it groups} with $\Gamma$-action. If $f$ is a unit in $L^1 (\Gamma)$, then $ff^*$ is a unit in $L^1 (\Gamma)$ as well and it is positive in $\Nh\Gamma$. Thus, under condition a) theorem \ref{t1_neu} gives
\[
h_{ff^*} = \log \ddet_{\Nh\Gamma} ff^* = 2 \log \ddet_{\Nh\Gamma} f \; .
\]
Combining this with the expected equation $2h_f = h_{ff^*}$ we would get the entropy formula $h_f = \log \ddet_{\Nh\Gamma} f$ under conditions a) and b) of theorem \ref{t1_neu} only. 

Note that in this argument we would need only a special case of the Yuzvinskii formula. In particular all actions are expansive by proposition \ref{t9}.
%\newpage
%\input{lit}

%\input{address}
Mathematisches Institut\\
Einsteinstr. 62\\
48149 M\"unster, Germany\\
deninger@math.uni-muenster.de

\begin{thebibliography}{9999}
\bibitem[B1]{B1}D. Boyd, Mahler's measure and special values of $L$-functions. Experiment. Math. {\bf 7} (1998), 37--82
\bibitem[B2]{B2}D. Boyd, Mahler's measure and invariants of hyperbolic manifolds. Number theory for the millenium I (Urbana IL), 127--143, A.K. Peters 2002
\bibitem[D]{D}C. Deninger, Deligne periods of mixed motives, $K$-theory and the entropy of certain $\Z^n$-actions. Journal of the AMS {\bf 10}, 259--281 (1997)
\bibitem[Di]{Di}J. Dixmier, $C^*$-algebras. North-Holland Math. Library Vol. 15, 1982
\bibitem[E]{E}G. Elek, Amenable groups, topological entropy and Betti numbers. Israel Journal of Math. {\bf 132} (2002), 315--335
\bibitem[FK]{FK}B. Fuglede, R.V. Kadison, Determinant theory in finite factors. Ann. Math. {\bf 55} (1952), 520--530
\bibitem[Gri]{Gri}R.I. Grigorchuk, On growth in group theory. Proc. ICM Kyoto 1990, Vol. I, 325--338
%\bibitem[Gro]{Gro}M. Gromov, Groups of polynomial growth and expanding maps. Publ IHES {\bf 53} (1981), 53--73
\bibitem[K]{K}Y. Katznelson, An introduction to harmonic analysis. Cambridge Univ. Press 2004
\bibitem[Le]{Le}H. Leptin, On one-sided harmonic analysis in noncommutative locally compact groups. J. Reine Angew. Math. {\bf 306} (1979), 122-153
\bibitem[LeP]{LeP}H. Leptin, D. Poguntke, Symmetry and nonsymmetry for locally compact groups. J. Funct. Anal. {\bf 33} (1979), 119-134
\bibitem[LSW]{LSW}D. Lind, K. Schmidt, T. Ward, Mahler measure and entropy for commuting automorphisms of compact groups. Invent. Math. {\bf 101} (1990), 593--629
\bibitem[LW]{LW1} E. Lindenstrauss, B. Weiss, Mean topological dimension. Israel Journal of Math. {\bf 115} (2000), 1--24
%\bibitem[LW2]{LW2}
\bibitem[Lo]{Lo}V. Losert, A characterization of groups with the one-sided Wiener property. J. Reine Angew. Math. {\bf 331} (1982), 47--57
\bibitem[L]{L}W. L\"uck, $L^2$-Invariants: Theory and applications to Geometry and $K$-Theory. Springer 2002
\bibitem[M]{M}J. Moulin Ollagnier, Ergodic theory and statistical mechanics. Springer LNM {\bf 1115}, 1985
\bibitem[OW]{OW}D. Ornstein, B. Weiss, Entropy and isomorphism theorems for actions of amenable groups. Journal d'Analyse Math\'ematique {\bf 48} (1987), 1--141
\bibitem[P]{P}A.L.T. Paterson, Amenability. AMS Math. Surveys and Monographs {\bf 29}, 1988
\bibitem[Sch]{Sch}K. Schmidt, Dynamical systems of algebraic origin. Birkh\"auser 1995
\bibitem[S]{S}R. Solomyak, On coincidence of entropies for two classes of dynamical systems. Ergod. Th. \& Dynam. Sys. {\bf 18} (1998), 731--738
\bibitem[W]{W}N. Wiener, Tauberian theorems. Ann. Math. {\bf 33} (1932), 1--100
\bibitem[Y]{Y}K. Yosida, Functional analysis. Springer 1971
\end{thebibliography}
\end{document}